\documentclass[12pt,twoside]{article}
\usepackage{amsmath,amssymb,amsthm}

\theoremstyle{plain}
\newtheorem{thm}{Theorem}[section]

\newtheorem{lem}{Lemma}[section]
\newtheorem{cor}{Corollary}[section]
\theoremstyle{remark}
\newtheorem{rem}{Remark}[section]

\newtheorem*{pf}{Proof}
\newtheorem*{prf1}{Proof of Theorem \ref{Theorem:W_as}}
\newtheorem*{prf2}{Proof of Theorem \ref{Theorem:fixv_degree_as}}
\newtheorem*{prf3}{Proof of Theorem \ref{Theorem:maxweight_as}}
\newtheorem*{prf4}{Proof of Theorem \ref{Theorem:maxdegree_as}}

\def\P{{\mathbb {P}}}                       
\def\E{{\mathbb {E}}}                       
                       
\def\N{{\mathbb {N}}}                       
\def\I{{\mathbb {I}}}

\def\FD{{\mathcal F}}

\def\o{{\rm {o}}} 
\def\O{{\rm {O}}} 

\def\CC{{\rm\kern.24em
   \vrule width.02em
       height1.4ex depth-.05ex
   \kern-.26em C}}
\def\QQ{{\rm\kern.24em
   \vrule width.02em
       height1.4ex depth-.05ex
   \kern-.26em Q}}
\def\PP{{\rm I\kern-.25em P}}         \def\RR{{\rm I\kern-.25em R}}
\def\DD{{\rm I\kern-.25em D}}         \def\EE{{\rm I\kern-.25em E}}
\def\FF{{\rm I\kern-.25em F}}         \def\NN{{\rm I\kern-.23em N}}
\def\RRp{{\rm I\kern-.25em R}_{+}}
\def\IND{{\rm 1\kern-.25em I}}

\begin{document}
%\thispagestyle{empty}
%\phantom{.}

\centerline{\large{\bf {The asymptotic behaviour of the weights and the degrees}}}
\smallskip
\centerline{\large{\bf {in an $N$-interactions random graph model}}}

\bigskip

\centerline{ {\sc Istv\'an Fazekas} and {\sc Bettina Porv\'azsnyik} }

\bigskip
Faculty of Informatics, University of Debrecen,  P.O. Box 12, 4010 Debrecen, Hungary, 
\centerline{
e-mail: fazekas.istvan@inf.unideb.hu, porvazsnyik.bettina@inf.unideb.hu.}

\medskip \bigskip

\begin{abstract}
A random graph evolution based on the interactions of $N$ vertices is studied.
During the evolution both the preferential attachment method and the uniform choice of vertices are allowed.
The weight of a vertex means the number of its interactions.
The asymptotic behaviour of the weight and the degree of a fixed vertex, moreover the limit of the maximal weight and the maximal degree are described.
The proofs are based on martingale methods.
\end{abstract}

\renewcommand{\thefootnote}{}
{\footnotetext{
{\bf Key words and phrases:}
Random graph, preferential attachment, scale-free, power law, submartingale, Doob-Meyer decomposition.

{\bf Mathematics Subject Classification:} 05C80, %Random graphs
60G42. %Martingales with discrete parameter 

Istv\'an Fazekas was supported by the T\'AMOP-4.2.2.C-11/1/KONV-2012-0001 project. The project has been supported by the European Union, co-financed by the European Social Fund.

Bettina Porv\'azsnyik was supported by the European Union and the State of Hungary, co-financed by the European Social Fund in the framework of T\'AMOP 4.2.4.A/2-11-1-2012-0001  \textsuperscript{`}National Excellence Program\textsuperscript{'}.
}

\section{Introduction}

Network theory is one of the most popular research topics.
During the last two decades, many types of networks were investigated.
Emprirical studies show that several real-world networks have certain similar features.
For an overview of random graph models and their properties see \cite{durrett} or \cite{hoff}.
It is known that many real life networks (the WWW, biological and social networks) are scale-free, see \cite{durrett}.
To describe the evolution of such networks, in \cite{barabasi} the preferential attachment model was suggested.

The scale-free property means that the asymptotic degree distribution follows a power law.
Besides the degree distribution, other characteristics are worth to study.
The degree of a fixed vertex and the maximal degree in preferential attachment models were investigated (see \cite{hoff}).
In \cite{Mo1}, the maximum degree in a general random tree model was examined, which includes the Barab\'asi-Albert random tree as a special case.
In \cite{Mo2} and \cite{Mo3},  the asymptotic behaviour of the degree of a given vertex and the maximal degree was studied in a $2$-parameter scale-free random graph model.
A well-known technique to analyze the growth of the maximal degree is M\'ori's martingale method (\cite{durrett}, \cite{Mo0}).

There are several modifications of the preferential attachment model (see \cite{hoff}, \cite{cooper}, \cite{SGWN}).
A random graph model based on the interactions of three vertices was introduced and power law degree distribution in that model was proved in \cite{BaMo1} and \cite{BaMo2}.
%The power law degree distribution in that model was proved in \cite{BaMo2}.
Instead of the three-interactions model, the interactions of $N$ vertices ($N \geq 3$ fixed) were studied in \cite{FIPB2}.
Scale-free properties for these generalized models were obtained in \cite{FIPB} and \cite{FIPB2}.

In this paper we extend some results of \cite{BaMo1} and \cite{BaMo2} to the $N$-interactions model.
Our aim is to study the asymptotic behaviour of the weight and the degree of a fixed vertex.
Moreover, we shall consider the limiting properties of the maximal weight and the maximal degree, as well.
In our proofs we follow the lines of \cite{BaMo1} and \cite{BaMo2}.
\\
\\
\textit{The $N$-interactions model}
\\
\\
In this paper we study the following $N$-interactions model (see \cite{FIPB2}).
A complete graph with $m$ vertices we call an $m$-clique, for short.
We denote an $m$-clique by the symbol $K_m$.
At time $n=0$ we start with a $K_N$.
The initial weight of this graph  and the initial weight all of its cliques are one. 
%This graph contains $N$ vertices, $ \binom{N}{2}$ edges, \dots , $\binom{N}{M}$ $M$-cliques $\left( M \leq N \right)$. 
%Each of these objects has initial weight $1$.
After the initial step we start to increase the size of the graph. 
At each step, the evolution of the graph is based on the interaction of $N$ vertices. 
More precisely, at each step $n=1,2,\dots$ we consider $N$ vertices and draw all non-existing edges between these vertices.
So we obtain a $K_N$.
The weight of this graph $K_N$ and the weights of all cliques in $K_N$ are increased by $1$. 
%(That is we increase the weights of $N$ vertices, $ \binom{N}{2}$ edges, \dots , $N$ different $\left(N-1\right)$-cliques and the $N$-clique $K_N$ itself.)
The choice of the $N$ interacting vertices is the following.

There are two possibilities at each step.
With probability $p$ we add a new vertex that interacts with $N-1$ old vertices, on the other hand,
with probability $\left( 1-p \right)$, $N$ old vertices interact. 
Here $0 < p \leq 1$ is fixed.

When we add a new vertex, then we choose $N-1$ old vertices and they together will form an $N$-clique.
However, to choose the $N-1$ old vertices we have two possibilities.
With probability $r$ we choose an $\left(N-1\right)$-clique from the existing $\left(N-1\right)$-cliques according to the weights of the $\left(N-1\right)$-cliques.
It means that an $\left(N-1\right)$-clique of weight $w_t$ is chosen with probability $w_t/\sum_h w_h$.
On the other hand, with probability $1-r$, we choose among the existing vertices uniformly, that is all $N-1$ vertices have the same chance.

At a step when we do not add a new vertex, then $N$ old vertices interact. 
As in the previous case, we have two possibilities. 
With probability $q$, we choose one $K_N$ of the existing $N$-cliques according to their weights.
It means that an $N$-clique of weight $w_t$ is chosen with probability $w_t/\sum_h w_h$.
On the other hand, with probability $1-q$, we choose among the existing vertices uniformly, that is all subsets consisting of $N$ vertices have the same chance.

In this paper we describe the asymptotic behaviour of the weight and the degree of a fixed vertex (Theorems \ref{Theorem:W_as} and \ref{Theorem:fixv_degree_as}), moreover we find the limit of the maximal weight and the maximal degree (Theorems \ref{Theorem:maxweight_as} and \ref{Theorem:maxdegree_as}).
The theorems are listed in Section $2$.
All the proofs and some important auxiliary results are presented in Section $3$.

\section{Main results}
Let us introduce the following notations.
$$
\alpha_1 = \left(1-p\right) q, \quad  \alpha_2 = \dfrac{N-1}{N}pr,  \quad  \alpha = \alpha_1 + \alpha_2,
$$
$$
\beta_1 =  \left( N-1 \right)\left( 1-r\right), \quad  \beta_2 = \dfrac{N\left( 1-p \right)\left( 1-q \right)}{p}, \quad
\beta =  \beta_1 + \beta_2.
$$

At time $n=0$ the initial complete graph on $N$ vertices is symmetric in the sense that all of its $N$ vertices have the same weight.
Let one of these $N$ vertices be labelled by $0$.
The other $N-1$ vertices are not labelled.
When a new vertex is born, let it be labelled by $1,2,\dots$ according the order in which they are added.
Let $j \geq 0$ be fixed integer.
Assume that the $j$th vertex exists after $l$ steps, where $0\leq j\leq l \leq n$.
Let us denote by $W[n,j]$ the weight of the $j$th vertex after the $n$th step.
Let us denote by $D[n,j]$ the degree of the $j$th vertex after $n$ steps.
(If $n < j$, then $W[n,j]=D[n,j]=0$.)
%(If $n < j$, then $W[n,j]=0$.)

\subsection{The weight and the degree of a fixed vertex}
The following theorem describes the asymptotic behaviour of the weight of a fixed vertex.
It is an extension of \textit{Theorem 4.1} in \cite{BaMo1}.
\begin{thm}   \label{Theorem:W_as}
Let $j \geq 0$ be fixed and let $\alpha >0$.
Then
\begin{equation} \label{W_as}
W[n,j] \sim  \dfrac{1}{\Gamma\left(1+\alpha\right)}\gamma_j n^{\alpha}
\end{equation}
almost surely as $n \to \infty$, where $\gamma_j$ is a positive random variable.
\end{thm}
We  turn to the asymptotic behaviour of the degree of a fixed vertex.
The following theorem is an extension of \textit{Theorem 5.2} in \cite{BaMo2}.
\begin{thm}   \label{Theorem:fixv_degree_as}
Let $j \geq 0$ be fixed and let $\alpha > 0$.
Then we have
\begin{equation} \label{D_as}
D[n,j] \sim  \dfrac{1}{\Gamma\left(1+\alpha\right)}\dfrac{\alpha_2}{\alpha}\gamma_j n^{\alpha}
\end{equation}
almost surely as $n \to \infty$, where the $\gamma_j$ positive random variable is defined in \eqref{W_as}.
\end{thm}
\subsection{The maximal weight and the maximal degree}
In this subsection we use the following notations. At time $n=0$, the vertices of the initial complete graph on $N$ vertices are labelled by $0,-1,\dots,-(N-1)$.
Let us denote by $\mathcal{W}_n$ the maximum of the weights of the vertices after $n$ steps, that is
\begin{equation} \label{maxweight_def}
\mathcal{W}_n = max\{W[n,j]:-(N-1) \leq j \leq n  \}.
\end{equation}
The following theorem is an extension of \textit{Theorem 5.1} in \cite{BaMo2} for the $N$-interactions model.
\begin{thm}   \label{Theorem:maxweight_as}
Let $\alpha > 0$. Then we have
\begin{equation} \label{maxweight_as}
\mathcal{W}_n \sim \dfrac{1}{\Gamma\left(1+\alpha\right)}\mu n^{\alpha} \text{ almost surely, as } n \to \infty,
\end{equation}
where $\mu = sup\{\gamma_j : j \geq -\left(N-1\right) \}$ is a finite positive random variable with $\gamma_j$ defined in Theorem \ref{Theorem:W_as}.
\end{thm}
Let us denote by $\mathcal{D}_n$ the maximal degree after $n$ steps, that is
\begin{equation} \label{maxdegree_def}
\mathcal{D}_n = max\{D[n,j]:-(N-1) \leq j \leq n  \}.
\end{equation}
The following theorem is an extension of \textit{Theorem 5.3} in \cite{BaMo2}.
\begin{thm}   \label{Theorem:maxdegree_as}
Let $\alpha > 0$. Then we have
\begin{equation} \label{maxdegree_as}
\mathcal{D}_n \sim \dfrac{1}{\Gamma\left(1+\alpha\right)}\dfrac{\alpha_2}{\alpha}\mu n^{\alpha} \text{ almost surely as } n \to \infty,
\end{equation}
where $\mu = sup\{\gamma_j : j \geq -\left(N-1\right) \}$ is the positive random variable defined in Theorem \ref{Theorem:maxweight_as}.
\end{thm}
\section{Proofs and auxiliary lemmas}
Introduce the following notations.
Let $\FD_{n-1}$ denote the $\sigma$-algebra of observable events after the $(n-1)$th step.
Let $j \geq 0$ be a fixed integer.
$W[n,j]$ is the weight of the $j$th vertex after the $n$th step.
(If $n < j$, then $W[n,j]=0$.)
Let $I[n,j]$ be the indicator of the event that the $j$th vertex exists after $n$ steps, that is 
$$
I[n,j]
=
\left\{
\begin{array}{ll}
1 ,
&
\textrm{if $W[n,j]>0$,}\\
0 ,
&
\textrm{if $W[n,j]=0$.}
\end{array}
\right.
$$
Let $J[n,j]$ be the indicator of the event that the $j$th vertex is born at the $n$th step.
Then $J[n,j]= I[n,j]-I[n-1,j]$.
For all $j,k,l,\, 0\leq j \leq l,\, 1 \leq k$, fixed positive integers, we consider the following sequences:
\begin{equation} \label{b_def}
b[n,k] = \prod_{\substack{i=1}}^n \left( 1 + \dfrac{\alpha k}{i} \right)^{-1},
\end{equation}
\begin{equation} \label{d_def}
d[n,k,j] = -\sum_{\substack{i=1}}^{n-1} b[i+1,k] \dfrac{\beta p}{V_i} \binom{W[i,j]+k-1}{k-1},
\end{equation}
\begin{equation}  \label{e_def}
e_n = \prod_{i=1}^n\left( 1-\dfrac{\alpha}{i}\, \right)^{-1}.
\end{equation}

Here we can see that the sequences $b[n,k]$ and $e_n$ are deterministic, while $d[n,k,j]$ is a sequence of $\FD_{n-1}$-measurable random variables for any $k$ and $j$.
Using the definition of $b[n,k]$ and the Stirling-formula for the Gamma function, we can show that
\begin{equation} \label{b_as}
b[n,k] \sim b_k n^{-k\alpha},\quad \text{ as } n \to \infty,
\end{equation}
where $b_k=\Gamma\left(1+\alpha k\right) > 0$, $k$ is fixed.
Moreover, we can easily see that
\begin{equation} \label{e_as}
e_n \sim \Gamma\left(1-\alpha\right) n^{\alpha},\quad \text{ as } n \to \infty.
\end{equation}

In the following lemma we introduce a martingale which will play important role in the proofs.
This lemma is an analouge of \textit{Lemma 4.1} in \cite{BaMo1} and \textit{Lemma 5.1} in \cite{BaMo2}.

\begin{lem}  \label{Lemma:Z_nkj}
Let $j,k,l, \, 0 \leq j \leq l$ be fixed nonnegative integers and let
\begin{equation} \label{Z_nkj}
Z[n,k,j] = \left(  b[n,k]\binom{W[n,j]+k-1}{k} + d[n,k,j] \right) I[l,j].
\end{equation}
Then $\left(Z[n,k,j], \FD_{n}\right)$ is a martingale for $n \geq l$.
\end{lem}
\begin{pf}
At each step, the weight of a fixed vertex is increased by $1$ if and only if it takes part in an interaction.
As in \cite{BaMo1}, \cite{FIPB} and in \cite{FIPB2}, it is easy to show that the conditional probability that the $j$th vertex takes part in an interaction at step $(n+1)$ is
\begin{equation}  \label{tpi}
\dfrac{W[n,j]}{n+1}\alpha + \dfrac{p}{V_{n}}\beta,
\end{equation}
provided that $W[n,j] > 0$.
Using this, we can see for $n \geq l$
\begin{equation*}
{\E} \left\{ \binom{W[n+1,j]+k-1}{k} I[l,j]| \FD_n \right\} =
\end{equation*}
\begin{equation*}
=
I[l,j] \left( 1- \left( \dfrac{W[n,j]}{n+1}\alpha + \dfrac{p}{V_{n}}\beta \right)  \right) \binom{W[n,j]+k-1}{k} +
I[l,j] \left( \dfrac{W[n,j]}{n+1}\alpha + \dfrac{p}{V_{n}}\beta \right) \binom{W[n,j]+k}{k} =
\end{equation*}
\begin{equation*}
=
I[l,j] \dfrac{p}{V_n}\beta\binom{W[n,j]+k-1}{k-1} +
I[l,j] \left( 1+\alpha \dfrac{k}{n+1}  \right)  \binom{W[n,j]+k-1}{k}.
\end{equation*}
Multiplying both sides by $b[n+1,k]$, we obtain
\begin{equation*}
{\E} \left\{ b[n+1,k]\binom{W[n+1,j]+k-1}{k} I[l,j]| \FD_n \right\} =
\end{equation*}
\begin{equation}    \label{W1}
=
I[l,j] \left( \binom{W[n,j]+k-1}{k} b[n,k] + d[n,k,j] - d[n+1,k,j] \right).
\end{equation}
Using that $d[n+1,k,j]$ is $\FD_n $-measurable, we obtain the desired result.
\hfill$\Box$
\end{pf}
\begin{lem} \label{Lemma:e_sm}
\begin{equation} \label{e_sm}
\left(  \dfrac{e_n I[k,j]}{W[n,j]-1}, \FD_n \right) \text{ is a nonnegative supermartingale, where }n=j,j+1,\dots \quad .
\end{equation}
\end{lem}
\begin{pf}
In a similar way as in the proof of Lemma \ref{Lemma:Z_nkj}, we have for $n \geq k$
\begin{equation*}
{\E} \left\{ \dfrac{I[k,j]}{W[n+1,j]-1}| \FD_n \right\} = 
\end{equation*}
\begin{equation*}
=
\left( \dfrac{W[n,j]}{n+1}\alpha + \dfrac{p}{V_{n}}\beta \right)  \dfrac{I[k,j]}{W[n,j]}  +
 \left( 1- \left( \dfrac{W[n,j]}{n+1}\alpha + \dfrac{p}{V_{n}}\beta \right)\right)\dfrac{I[k,j]}{W[n,j]-1}
=
\end{equation*}
\begin{equation*}
=
\left( \dfrac{W[n,j]}{n+1}\alpha + \dfrac{p}{V_{n}}\beta \right) 
\left( \dfrac{I[k,j]}{W[n,j]} - \dfrac{I[k,j]}{W[n,j]-1} \right) +
\dfrac{I[k,j]}{W[n,j]-1} \leq
\end{equation*}
\begin{equation} \label{fkor}
\leq
-\dfrac{\alpha I[k,j]}{(n+1)(W[n,j]-1)} + \dfrac{I[k,j]}{W[n,j]-1}
=
\dfrac{I[k,j]}{W[n,j]-1} \left( 1-\dfrac{\alpha}{n+1}  \right).
\end{equation}
Here we used that
$$
\left( \dfrac{W[n,j]}{n+1}\alpha + \dfrac{p}{V_{n}}\beta \right)
I[k,j]
\left(- \dfrac{1}{W[n,j](W[n,j]-1)} \right) \leq 
-\dfrac{\alpha I[k,j]}{(n+1)(W[n,j]-1)}.
$$
Multiplying both sides of \eqref{fkor} by $e_{n+1}$, we obtain the result.
\hfill$\Box$
\end{pf}
\begin{prf1}
Following the method presented in \cite{BaMo1}, the proof contains two parts.
First, we will show that the above result is valid with nonnegative $\gamma_j$.
In the second part of the proof we will show that $\gamma_j$ is positive with probability $1$.

Let $B_{n+1} = \lbrace{ W[n+1,j] = W[n,j]+1  \rbrace} $.
Consider the event that the $j$th vertex exists after $n$ steps.
On this event, by \eqref{tpi},
\begin{equation} \label{akor}
{\P} \left( B_{n+1} | \FD_n \right) \geq \dfrac{\alpha}{n+1} .
\end{equation}
The sequence $\left( B_n, n \in {\N}\right)$ is adapted to the sequence $\left(\FD_n, n \in {\N}\right)$ of $\sigma$-algebras.
Using \textit{Corollary VII-2-6} of \cite{neveu} and \eqref{akor}, we have
\begin{equation} \label{W_infty}
W[n,j] \to \infty \quad \text{ a.s. as } n \to \infty.
\end{equation}
Consider the martingale $\left(Z[n,k,j], \FD_{n}\right)$ in Lemma \ref{Lemma:Z_nkj} and let $k=1$.
Then
\begin{equation}   \label{Z1}
Z[n,1,j] = \left( b[n,1] W[n,j] + d[n,1,j]  \right) I[l,j].
\end{equation}
Applying the Marcinkiewicz strong law of large numbers to the number of vertices, we have
\begin{equation}  \label{Markinkiewicz}
V_n = pn + \o \left( n^{1/2 + \varepsilon} \right)
\end{equation}
almost surely, for any $\varepsilon >0$.
By this and \eqref{b_as}, we obtain that
\begin{equation*}
d[n,1,j] =  -\sum_{\substack{i=1}}^{n-1} b[i+1,1] \dfrac{\beta p}{V_i} \sim
 -\beta \Gamma\left(1+\alpha\right) \sum_{\substack{i=1}}^{n-1} \dfrac{p}{pi + \o \left( i^{1/2 + \varepsilon}\right) }i^{-\alpha} =
\end{equation*}
\begin{equation*}
=
-\beta \Gamma\left(1+\alpha\right) \sum_{\substack{i=1}}^{n-1} i^{-\left(1+\alpha\right)} \left(1+\o \left( 1 \right) \right).
\end{equation*}
Using that $\alpha>0$, we see that $d[n,1,j]$ converges, as $n \to \infty$, and therefore the martingale $Z[n,1,j]$ is bounded from below.
Moreover, we shall see that the martingale $Z[n,1,j]$ has bounded differences.
The sequence $b[n,1]$ is monotonically decreasing, hence
\begin{equation*}
Z[n+1,1,j]-Z[n,1,j] \leq b[n,1] \left(W[n+1,j] - W[n,j] \right) \leq 1.
\end{equation*}
It is also easy to compute that
\begin{equation*}
Z[n,1,j]-Z[n+1,1,j] \leq \left(b[n,1] - b[n+1,1]\right) W[n,j] + \left(d[n,1,j] -d[n+1,1,j]\right) =
\end{equation*}
\begin{equation*}
=
\left(b[n,1] -b[n+1,1]\right) W[n,j] + b[n+1,1]\dfrac{p}{V_{n}}\beta \leq
b[n+1,1]\left(\alpha + \dfrac{p}{N}\beta  \right) \leq \alpha + \dfrac{p}{N}\beta .
\end{equation*}
As the martingale $Z[n,1,j]$ is bounded from below and it has bounded differences, by \textit{Proposition VII-3-9} of \cite{neveu}, it is convergent almost surely, as $n \to \infty$.
By the definition of $Z[n,1,j]$, we see that $b[n,1]W[n,j]$ also converges almost surely on the event $\{ W[l,j] > 0  \}$.
This, \eqref{W_infty}, and \eqref{b_as} implies that \eqref{W_as} is true with nonnegative $\gamma_j$.

Now we will show that $\gamma_j$ is positive with probability $1$.
Consider the supermartingale $\label{e_sm}
\left(  \dfrac{e_n I[k,j]}{W[n,j]-1}, \FD_n \right), n\geq j$, in Lemma \ref{Lemma:e_sm}.
This supermartingale is nonnegative, therefore according to the submartingale-convergence theorem, it converges almost surely.
\[ \lim_{l \to \infty} I[l,j] = 1\] almost surely, hence $\dfrac{e_n }{W[n,j]-1}$ also converges almost surely, as $n \to \infty$.
This and \eqref{e_as} imply that $\gamma_j$ is positive almost surely.
\hfill$\Box$
\end{prf1}
Assume that the $j$th vertex exists after $l$ steps, where $0\leq j\leq l \leq n$.
We denote by $D[n,j]$ the degree of the $j$th vertex after $n$ steps ($j \geq 0$ is fixed).
\begin{rem} \label{Remark}
Contrary to the weights, at each step the degree of a fixed vertex can grow by $0,1,\dots, N-1$.
Hence $0 \leq D[n,j]-D[n-1,j] \leq N-1$ for all fixed $j \geq 0$.
Moreover, the degree of a fixed vertex does not change in such interactions when we do not add a new vertex and the choice is done according to the preferential attachment rule.
\end{rem}
\begin{lem} \label{Lemma:CondExp}
Let $j\geq 0$ be fixed. We have
\begin{equation*}
{\E} \{ I[k,j]D[n+1,j]|\FD_{n} \} =
I[k,j]\left(
D[n,j] \left[ \left(1-\dfrac{W[n,j]}{n+1}\alpha - \dfrac{p}{V_{n}}\beta \right)+ \right.\right.
\end{equation*}
\begin{equation*}
+ \left.
 \left( 1-p \right) \left( q\dfrac{W[n,j]}{n+1} + \left( 1-q \right)\dfrac {\binom{D[n,j]}{N-1}} {\binom{V_{n}}{N}} \right) \right]  +
\end{equation*}
\begin{equation*}
+ 
\left(D[n,j]+1\right)  \left[ p \left(r \dfrac{\left(N-1\right) W[n,j]}{N\left(n+1\right)} + \left( 1-r \right) \dfrac{\binom{D[n,j]}{N-2}}{\binom{V_{n}}{N-1}}  \right) + \right.
\end{equation*}
\begin{equation*}
\left.
+ \left( 1-p \right)\left( 1-q \right)\dfrac {\binom{D[n,j]}{N-2} \left( V_{n}-D[n,j]-1 \right)} {\binom{V_{n}}{N}} \right] + \dots +
\end{equation*}
\begin{equation*}
+ 
\left(D[n,j]+m\right)\left[ p \left( 1-r \right)\dfrac { \binom{ 
D[n,j]}{N-m-1} \binom{V_{n}-D[n,j]-1}{m-1}} {\binom{V_{n}}{N-1}} + \right.
\end{equation*}
\begin{equation*}
\left.
+ \left( 1-p \right)\left( 1-q \right)\dfrac {\binom{D[n,j]}{N-m-1} \binom{ V_{n}-D[n,j]-1}{m} } {\binom{V_{n}}{N}} \right] + \dots +
\end{equation*}
\begin{equation*}
+ 
\left(D[n,j]+(N-1)\right) \left[ p\left( 1-r \right) \dfrac{\binom{V_{n}-D[n,j]-1}{N-2}}{\binom{V_{n}}{N-1}} \right.
\end{equation*}
\begin{equation} \label{CondExp}
\left.
\left. +
\left( 1-p \right) \left( 1-q \right)\dfrac{\binom{V_{n}-D[n,j]-1}{N-1}}{\binom{V_{n}}{N}} \right]\right),
\end{equation}
for $n \geq k$.
\end{lem}
\begin{pf}
As in \cite{FIPB2}, the probability that an old vertex of weight $W[n,j]$ takes part in the interaction at step $(n+1)$ is
$$
\dfrac{W[n,j]}{n+1}\alpha + \dfrac{p}{V_{n}}\beta. 
$$

Consider a fixed vertex with weight $W[n,j]$ and degree $D[n,j]$. 
Using the basic properties of the model, we have the probability that in the $(n+1)$th step 
\begin{itemize}
 \item its degree $D[n,j]$ does not change is
 $$
1-\left(\dfrac{W[n,j]}{n+1}\alpha + \dfrac{p}{V_{n}}\beta \right)+\left( 1-p \right) \left( q\dfrac{W[n,j]}{n+1} + \left( 1-q \right)\dfrac {\binom{D[n,j]}{N-1}} {\binom{V_{n}}{N}} \right)\,;
 $$
 \item  its degree is increased by 1 is
 $$
p \left(r \dfrac{\left(N-1\right)W[n,j]}{N(n+1)} + \left( 1-r \right) \dfrac{\binom{D[n,j]}{N-2}}{\binom{V_{n}}{N-1}}  \right) +
\left( 1-p \right)\left( 1-q \right)\dfrac {\binom{D[n,j]}{N-2} \left( V_{n}-D[n,j]-1 \right)} {\binom{V_{n}}{N}}\,;
 $$
 \item its degree is increased by $m$ ($1<m \leq N-1$) is
 $$
p \left( 1-r \right)\dfrac { \binom{D[n,j]}{N-m-1} \binom{V_{n}-D[n,j]-1}{m-1}} {\binom{V_{n}}{N-1}} + \left( 1-p \right)\left( 1-q \right)\dfrac {\binom{D[n,j]}{N-m-1} \binom{ V_{n}-D[n,j]-1}{m} } {\binom{V_{n}}{N}}\,.
 $$
% \item its degree is increased by $N-1$ is
% $$
%p \left( 1-r \right) \dfrac{\binom{V_{n}-D[n,j]-1}{N-2}}{\binom{V_{n}}{N-1}} +
%\left( 1-p \right) \left( 1-q \right)\dfrac{\binom{V_{n}-D[n,j]-1}{N-1}}{\binom{V_{n}}{N}}\,.
% $$
\end{itemize}
Using the above formulae, we obtain equation  \eqref{CondExp}.
\hfill$\Box$
\end{pf}
\begin{cor}  \label{Cor}
Let $j\geq 0$ be fixed. For $n \geq k$, we have
\begin{equation}
0 \leq {\E} \{ I[k,j]D[n+1,j]|\FD_{n} \} = I[k,j]\left( D[n,j] + \alpha_2 \dfrac{W[n,j]}{n+1} + R_n\right),
\end{equation}
where $0 \leq R_n \leq (N-1)\dfrac{p\beta}{V_n}$.
\end{cor}
\begin{pf}
Let $1 < m < N$ be an integer.
Using the above notations we can rewrite \eqref{CondExp} into the following form:
\begin{equation*}
{\E} \{ I[k,j]D[n+1,j]|\FD_{n} \} =
I[k,j] \left[ D[n,j] + \alpha_2 \dfrac{W[n,j]}{n+1} + R_n \right],
\end{equation*}
where
\begin{equation*}
R_n =
\beta_1 p \left[-\dfrac{D[n,j]}{V_n} + \left(D[n,j] +1 \right) \dfrac{1}{N-1} \dfrac{\binom{D[n,j]}{N-2}}{\binom{V_{n}}{N-1}} + \dots +
  \right.
\end{equation*}
\begin{equation*}
\left.
\dots +\left(D[n,j] +m \right) \dfrac{1}{N-1} \dfrac { \binom{D[n,j]}{N-m-1} \binom{V_{n}-D[n,j]-1}{m-1}} {\binom{V_{n}}{N-1}}  + \dots +
  \right.
\end{equation*}
\begin{equation*}
\left.
+ \dots +\left(D[n,j] +(N-1) \right) \dfrac{1}{N-1} \dfrac{\binom{V_{n}-D[n,j]-1}{N-2}}{\binom{V_{n}}{N-1}}
  \right] +
\end{equation*}
\begin{equation*}
+
\beta_2 p \left[-\dfrac{D[n,j]}{V_n} + \left(D[n,j]\right) \dfrac{1}{N} \dfrac{\binom{D[n,j]}{N-1}}{\binom{V_{n}}{N}}
+ \left(D[n,j] +1 \right) \dfrac{1}{N} \dfrac {\binom{D[n,j]}{N-2} \left( V_{n}-D[n,j]-1 \right)} {\binom{V_{n}}{N}} + \dots +
  \right.
\end{equation*}
\begin{equation*}
\left.
\dots +\left(D[n,j] +m \right) \dfrac{1}{N} \dfrac {\binom{D[n,j]}{N-m-1} \binom{ V_{n}-D[n,j]-1}{m} } {\binom{V_{n}}{N}}  + \dots +
  \right.
\end{equation*}
\begin{equation*}
\left.
+ \dots +\left(D[n,j] +(N-1) \right) \dfrac{1}{N} \dfrac{\binom{V_{n}-D[n,j]-1}{N-1}}{\binom{V_{n}}{N}}
  \right] = 
\end{equation*}
\begin{equation}   \label{R}
= \beta_1 p R^{(1)} + \beta_2 p R^{(2)}.
\end{equation}
Now, we give upper bounds for $R^{(1)}$ and $R^{(2)}$ separately.
It is easy to see that

\begin{equation*}
 R^{(1)}V_n = \left[-D[n,j] + \left(D[n,j] +1 \right) \dfrac{\binom{D[n,j]}{N-2}}{\binom{V_{n}-1}{N-2}} + \dots +
  \right.
\end{equation*}
\begin{equation*}
\left.
\dots +\left(D[n,j] +m \right) \dfrac { \binom{D[n,j]}{N-m-1} \binom{V_{n}-D[n,j]-1}{m-1}} {\binom{V_{n}-1}{N-2}}  + \dots +
  \right.
\end{equation*}
\begin{equation*}
\left.
+ \dots +\left(D[n,j] +(N-1) \right) \dfrac{\binom{V_{n}-D[n,j]-1}{N-2}}{\binom{V_{n}-1}{N-2}}
  \right] = 
\end{equation*}
\begin{equation*}
  = \left[-D[n,j] + D[n,j] \dfrac{1}{\binom{V_{n}-1}{N-2}} \sum_{k=0}^{N-2}
 \binom{D[n,j]}{k} \binom{V_{n}-D[n,j]-1}{N-2-k} +  \right.
\end{equation*}
\begin{equation*}
 + \left. \dfrac{1}{\binom{V_{n}-1}{N-2}} \sum_{k=0}^{N-2} \left( N-1-k \right)
 \binom{D[n,j]}{k} \binom{V_{n}-D[n,j]-1}{N-2-k} \right] \leq
\end{equation*}
\begin{equation} \label{R_1fkor}
 \leq -D[n,j] + D[n,j] + \left(N-1\right) =  \left(N-1\right) .
\end{equation}
Similarly, we have
\begin{equation*}
R^{(2)}V_n = \left[-D[n,j] + \left(D[n,j]\right)  \dfrac{\binom{D[n,j]}{N-1}}{\binom{V_{n}-1}{N-1}}
+ \left(D[n,j] +1 \right) \dfrac {\binom{D[n,j]}{N-2} \left( V_{n}-D[n,j]-1 \right)} {\binom{V_{n}-1}{N-1}} + \dots +
  \right.
\end{equation*}
\begin{equation*}
\left.
\dots +\left(D[n,j] +m \right) \dfrac {\binom{D[n,j]}{N-m-1} \binom{ V_{n}-D[n,j]-1}{m} } {\binom{V_{n}-1}{N-1}}  + \dots +
  \right.
\end{equation*}
\begin{equation*}
\left.
+ \dots +\left(D[n,j] +(N-1) \right) \dfrac{\binom{V_{n}-D[n,j]-1}{N-1}}{\binom{V_{n}-1}{N-1}}
  \right] =
\end{equation*}
\begin{equation*}
  = \left[-D[n,j] + D[n,j] \dfrac{1}{\binom{V_{n}-1}{N-1}} \sum_{k=0}^{N-1}
 \binom{D[n,j]}{k} \binom{V_{n}-D[n,j]-1}{N-1-k} +  \right.
\end{equation*}
\begin{equation*}
 + \left. \dfrac{1}{\binom{V_{n}-1}{N-1}} \sum_{k=0}^{N-1} \left( N-1-k \right)
 \binom{D[n,j]}{k} \binom{V_{n}-D[n,j]-1}{N-1-k} \right] \leq
\end{equation*}
\begin{equation} \label{R_2fkor}
 \leq -D[n,j] + D[n,j] + \left(N-1\right) =  \left(N-1\right) .
\end{equation}
Using \eqref{R}, \eqref{R_1fkor} and \eqref{R_2fkor}, we have 
\begin{equation}
R_n = \beta_1 p R^{(1)} + \beta_2 p R^{(2)} \leq \left(N-1\right)\dfrac{p\beta }{V_n}.
\end{equation}
The proof is complete.
\hfill$\Box$
\end{pf}
\begin{prf2}
Consider the following bounded random variable: $\xi_n= \dfrac{I[k,j]}{N-1} (D[n,j]-D[n-1,j])$.
By Remark \ref{Remark}, we have $0 \leq \xi_n \leq 1$.
Applying an appropriate version of \textit{Corollary VII-2-6} of \cite{neveu} (see \textit{Proposition 2.4} of \cite{backhausz}), then using Corollary \ref{Cor} and \eqref{W_as}, we have
\begin{equation*}
D[n,j] = (N-1)\sum_{i=1}^{n}\xi_i \sim (N-1)\sum_{i=1}^{n} {\E} \left(\xi_i | \FD_{i-1}\right) = \sum_{i=1}^{n} \left( \alpha_2 \dfrac{W[i-1,j]}{i} +R_{i-1} \right) \sim
\end{equation*}
\begin{equation}
\sim \dfrac{1}{\Gamma\left(1+\alpha\right)} \dfrac{\alpha_2}{\alpha} \gamma_j n^{\alpha},
\end{equation}
provided that $j$th vertex exists after $k$ steps.
As $\lim_{k \to \infty} W[k,j] = \infty$ a.s., we obtain the statement.
\hfill$\Box$
\end{prf2}
The following lemma is an extension of \textit{Lemma 5.2} in \cite{BaMo2}.
This statement is useful when we consider the asymptotic behaviour of the maximal weight.
\begin{lem}     \label{Lemma:S_mnk}
For all $k \geq 0$, $1 \leq m \leq n$ fixed nonnegative integers, let 
\begin{equation}   \label{S_mnk_def}
S[m,n,k] = \sum_{j=m}^{n} {\E} \left( b[n,k]  \binom{W[n,j]+k-1}{k} I[n,j] \right).
\end{equation}
Then
\begin{equation}  \label{S_mnk_korlat}
S[m,n,k] \leq C_k \sum_{j=m}^{n} j^{-\alpha k}
\end{equation}
with a positive constant $C_k$.
\end{lem}
\begin{pf}
As in \cite{BaMo2}, we use induction on $k$.
Let $k=0$.
Then
\begin{equation*}
S[m,n,0] = \sum_{j=m}^{n} {\E} \left( b[n,0] I[n,j] \right) = \sum_{j=m}^{n} {\P} \left( W[n,j] > 0 \right) \leq n-m+1\,.
\end{equation*}
Suppose that the statement is true for $k-1$, that is
\begin{equation}   \label{S_ihypothesis}
S[m,n,k-1] = \sum_{j=m}^{n} {\E} \left( b[n,k-1]  \binom{W[n,j]+k-2}{k-1} I[n,j] \right) \leq C_{k-1} \sum_{j=m}^{n} j^{-\alpha \left(k-1\right)}.
\end{equation}
By Lemma \ref{Lemma:Z_nkj}, $Z[n,k,j]$ is a martingale.
The difference of two martingales is also a martingale.
So, in the definition of $Z[n,k,j]$ changing $I[l,j]$ for $J[l,j]$, we obtain again a martingale.
Using the definitions of $J[n,j]$ and $Z[n,k,j]$, we have
\begin{equation*}
S[m,n,k] = \sum_{j=m}^{n} {\E} \left( \sum_{l=j}^{n}\left( b[n,k]  \binom{W[n,j]+k-1}{k} J[l,j] \right) \right)=
\end{equation*}
\begin{equation*}
 = 
 \sum_{j=m}^{n} {\E} \left( \sum_{l=j}^{n}\left( Z[n,k,j]-d[n,k,j]\right) J[l,j] \right) =
 \sum_{j=m}^{n} {\E} \left( \sum_{l=j}^{n}\left( Z[l,k,j]-d[n,k,j]\right) J[l,j] \right) =
\end{equation*}
\begin{equation} \label{2terms}
 = 
 {\E} \left( \sum_{j=m}^{n} \sum_{l=j}^{n} b[l,k] J[l,j] \right) +  {\E} \left( \sum_{j=m}^{n} \sum_{l=j}^{n} \left( d[l,k,j]-d[n,k,j]  \right) J[l,j] \right). 
\end{equation}
In the last step we used that $W(l,j)=1$ if $J[l,j]=1$.

Now, we give an upper bound for the two terms in \eqref{2terms} separately.
We have already seen that, for a fixed $k$, the sequence $b[n,k]$ is monotonically decreasing.
Therefore, applying also \eqref{b_as},
$$
{\E} \left( \sum_{j=m}^{n} \sum_{l=j}^{n} b[l,k] J[l,j] \right) \leq \sum_{j=m}^{n} b[j,k] {\E} \left(  \sum_{l=j}^{n} J[l,j] \right) \leq  \sum_{j=m}^{n} b[j,k] {\E} I[n,j]
\leq
$$
\begin{equation} \label{S1_k}
\leq
\sum_{j=m}^{n} b[j,k] \leq C_k^{(1)} \sum_{j=m}^{n} j^{-\alpha k}.
\end{equation}
For the second term in \eqref{2terms}, changing the order of summation and using that $I[i,j]=\sum_{l=j}^{i}J[l,j]$, we have
$$
{\E} \left( \sum_{j=m}^{n} \sum_{l=j}^{n} \left( d[l,k,j]-d[n,k,j]  \right) J[l,j] \right) =
{\E} \left( \sum_{j=m}^{n} \sum_{l=j}^{n}  \sum_{i=l}^{n-1}  b[i+1,k] \dfrac{p}{V_i}\beta \binom{W[i,j]+k-1}{k-1} J[l,j]     \right)=
$$
$$
=
{\E} \left( \sum_{i=m}^{n-1} \dfrac{b[i+1,k]}{b[i,k-1]}  \dfrac{p}{V_i}\beta \sum_{j=m}^{i}  b[i,k-1] \binom{W[i,j]+k-2}{k-1} \dfrac{W[i,j]+k-1}{W[i,j]} I[i,j]     \right) \leq
$$
\begin{equation}   \label{*}
\leq
k \sum_{i=m}^{n-1} \dfrac{b[i+1,k]}{b[i,k-1]} {\E} \left( \dfrac{p}{V_i}\beta  \sum_{j=m}^{i}  b[i,k-1] \binom{W[i,j]+k-2}{k-1}  I[i,j]     \right).
\end{equation}
In the last step we applied that $\dfrac{W[i,j]+k-1}{W[i,j]} \leq k$ if $I[i,j]>0$. 

Now, following the line of the proof presented in \cite{BaMo2}, we give upper bound on the events $\left\{ V_i < \dfrac{pi}{2} \right\}$ and $\left\{  V_i \geq \dfrac{pi}{2} \right\}$ separately.
Using the induction hypothesis \eqref{S_ihypothesis}, we have
$$
 {\E}  \left( \dfrac{p}{V_i}\beta {\I}_{\left\{  V_i \geq \frac{pi}{2} \right\}} \sum_{j=m}^{i}  b[i,k-1] \binom{W[i,j]+k-2}{k-1}  I[i,j]     \right)
 \leq
 \dfrac{2\beta}{i}S[m,i,k-1] \leq
$$
\begin{equation} \label{V_i>}
\leq \dfrac{2\beta}{i} C_{k-1} \sum_{j=m}^{i} j^{-\alpha \left(k-1\right)}.
\end{equation}
(Here ${\I}_A$ is the indicator of the set $A$.)
On the other hand, by \eqref{b_as},
$$
 {\E}  \left( \dfrac{p}{V_i}\beta {\I}_{\left\{ V_i < \frac{pi}{2} \right\}} \sum_{j=m}^{i}  b[i,k-1] \binom{W[i,j]+k-2}{k-1}  I[i,j]     \right)
 \leq
 $$
\begin{equation}  \label{V_i<}
\leq
\dfrac{p}{N}\beta  {\P} \left\{ V_i < \dfrac{pi}{2} \right\}\sum_{j=m}^{i}  b[i,k-1] \binom{i+k-2}{k-1} =
\O \left( e^{-\varepsilon i} i^{-(k-1)\alpha} i^{k-1}i \right) =
\o \left( \dfrac{1}{i} \sum_{j=m}^{i} j^{-\alpha \left(k-1\right)} \right),
\end{equation}
as $i \to \infty$.
In the above computation we used Hoeffding's exponential inequality (\textit{Theorem 2} in \cite{HO}) to obtain the following upper bound: ${\P} \left\{ V_i < \dfrac{pi}{2} \right\} \leq e^{-\varepsilon i}$, where $\varepsilon$ only depends on $p$, (actually $\varepsilon = \frac{p^2}{2}$).
Therefore, by \eqref{*}, \eqref{V_i<} and \eqref{V_i>}, we have
$$
{\E} \left( \sum_{j=m}^{n} \sum_{l=j}^{n} \left( d[l,k,j]-d[n,k,j]  \right) J[l,j] \right) \leq
$$
$$
\leq
k \sum_{i=m}^{n-1} \dfrac{b[i+1,k]}{b[i,k-1]} {\E} \left( \dfrac{p}{V_i}\beta  \sum_{j=m}^{i}  b[i,k-1] \binom{W[i,j]+k-2}{k-1}  I[i,j]     \right)
\leq
$$
\begin{equation}   \label{S2_k}
\leq
k C_k^{(2)'}\sum_{i=m}^{n-1} i^{-\alpha} \sum_{j=m}^{i} \dfrac{1}{i}j^{-\alpha \left(k-1\right)} 
\leq C_k^{(2)}\sum_{j=m}^{n} j^{-\alpha \left(k-1\right)} \sum_{i=j}^{n} i^{-\left(1+\alpha\right)} \leq 
C_k^{(3)}\sum_{j=m}^{n} j^{-\alpha k}.
\end{equation}
Above we applied that, by \eqref{b_as}, $\dfrac{b[i+1,k]}{b[i,k-1]} =  \O \left( i^{-\alpha} \right)$ as $i \to \infty$.
Finally, by \eqref{2terms}, \eqref{S1_k} and \eqref{S2_k}, we have
$$
S[m,n,k] \leq C_k^{(1)} \sum_{j=m}^{n} j^{-\alpha k} + C_k^{(3)}\sum_{j=m}^{n} j^{-\alpha k}.
$$
The proof is complete.
\hfill$\Box$
\end{pf}
\begin{prf3}
To obtain \eqref{maxweight_as}, we can apply the method \textit{Theorem 5.1} in \cite{BaMo2}.
Let
\begin{equation} \label{M[m,n]_def}
M[m,n] =max\{W[n,j]:-(N-1) \leq j < m  \},
\end{equation}
where $1 \leq m \leq n$ fixed.
From Theorem \ref{Theorem:W_as}, we have
\begin{equation}  \label{M}
\Gamma\left(1+\alpha\right)  \lim_{n \to \infty}  n^{-\alpha}M[m,n] = max\{\gamma_j:-(N-1) \leq j < m \}
\end{equation}
almost surely.
Using \eqref{W1}, we can prove that the following process is a submartingale:
\begin{equation*}
b[n,k]\binom{W[n,j]+k-1}{k} = b[n,k]\binom{W[n,j]+k-1}{k}I[n,j],\quad n\geq j.
\end{equation*}
Let $Q[m,n] = max_{m \leq j \leq n} W[n,j]$.
Then $0 \leq \mathcal{W}_n-M[m,n] \leq Q[m,n]$.
As the maximum of increasing numbers of submartingales is also a submartingale, therefore
\begin{equation*}
b[n,k]\binom{Q[m,n]+k-1}{k}, \quad n\geq m
\end{equation*}
is also a submartingale.
For nonnegative numbers the maximum is majorized by the sum.
Therefore, and by Lemma \ref{Lemma:S_mnk}, we obtain
\begin{equation}   \label{A}
{\E} \left( b[n,k]\binom{Q[m,n]+k-1}{k}\right) \leq S[m,n,k]\leq C_k\sum_{j=m}^{n} j^{-\alpha k}.
\end{equation}
Since
\begin{equation}  \label{B}
0 \leq \left(b[n,1]Q[m,n]\right)^k \leq 
\dfrac{b[n,1]^k}{b[n,k]}k!b[n,k]\binom{Q[m,n]+k-1}{k},
\end{equation}
we see that the submartingale $b[n,1] Q[m,n] $ is bounded in $L^k$ for all $k\alpha >1$.
Hence, this submartingale converges almost surely and in $L^k$ for every $k > \dfrac{1}{\alpha}$.
Moreover, by \eqref{b_as}, \eqref{A} and \eqref{B}, we have
\begin{equation}   \label{finite}
{\E} \left(\limsup_{n \to \infty}n^{-\alpha} Q[m,n]  \right)^k  \leq k! C_k \dfrac{1}{\Gamma\left(1+\alpha k\right)}\sum_{j=m}^{\infty} j^{-\alpha k}.
\end{equation}
Now, using the monoton convergence theorem, we have
\begin{equation*}
{\E} \left(\lim_{m \to \infty} \limsup_{n \to \infty}\left(n^{-\alpha}  Q[m,n]  \right)^k \right) = 0,
\end{equation*}
for $k > \frac{1}{\alpha}$.
As $Q[m,n]$ is decreasing, as $m$ increases, so
\begin{equation}   \label{Q}
\lim_{m \to \infty} \limsup_{n \to \infty} n^{-\alpha} Q[m,n] = 0\quad \text{ a.s.}
\end{equation}
Therefore, as $0 \leq \mathcal{W}_n-M[m,n] \leq Q[m,n]$,
\[
\lim_{m \to \infty} \limsup_{n \to \infty}\left(n^{-\alpha} \left( \mathcal{W}_n-M[m,n]  \right)\right) = 0\quad \text{ a.s.}
\]
This relation and \eqref{M} imply \eqref{maxweight_as}.
Using relation \eqref{finite} we can show that $\mu = sup\{\gamma_j : j \geq -\left(N-1\right) \}$ is a.s. finite.
\hfill$\Box$
\end{prf3}
\begin{prf4}
We follow the line of the proof of \textit{Theorem 5.3} in \cite{BaMo2}.
The evolution mechanism of the graph implies that $D[n,j] \leq (N-1)W[n,j]$.
Therefore we have
\begin{equation*}
max\{D[n,j]:-(N-1) \leq j < m  \} \leq \mathcal{D}_n \leq 
\end{equation*}
\begin{equation*}
 \leq max\{D[n,j]:-(N-1) \leq j < m  \} + max\{(N-1)W[n,j]:m \leq j \leq n  \}.
\end{equation*}
Multiplying both sides by $n^{-\alpha}$ and then considering the limit as $n \to \infty$, Theorem \ref{Theorem:fixv_degree_as} implies
\begin{equation*}
max\bigg\{\dfrac{1}{\Gamma\left(1+\alpha\right)}\dfrac{\alpha_2}{\alpha} \gamma_j:-(N-1) \leq j < m  \bigg\} \leq 
\end{equation*}
\begin{equation*}
\leq
\liminf_{n\to \infty} \mathcal{D}_n n^{-\alpha} \leq \limsup_{n\to \infty} \mathcal{D}_n n^{-\alpha} \leq 
\end{equation*}
\begin{equation*}
 \leq max\bigg\{\dfrac{1}{\Gamma\left(1+\alpha\right)}\dfrac{\alpha_2}{\alpha} \gamma_j:-(N-1) \leq j < m  \bigg\} + (N-1)\limsup_{n \to \infty} n^{-\alpha} Q[m,n]
\end{equation*}
as $n \to \infty$.
As $m \to \infty$, by \eqref{Q}, we obtain the desired result.
\hfill$\Box$
\end{prf4}
%%%

\end{document}